\documentclass[reqno]{amsart} 

\newtheorem{theorem}{Theorem}[section]
\newtheorem{proposition}{Proposition}[section]

\newtheorem{lemma}[theorem]{Lemma}
\theoremstyle{remark}
\newtheorem{remark}[theorem]{Remark}
\numberwithin{equation}{section}

\allowdisplaybreaks

\author{Michael Schlosser}
\address{Fakult\"at f\"ur Mathematik, Universit\"at Wien,
Nordbergstra{\ss}e 15, A-1090 Vienna, Austria}
\email{michael.schlosser@univie.ac.at}
\urladdr{http://www.mat.univie.ac.at/{\textasciitilde}schlosse}
\thanks{Partly supported by FWF Austrian Science Fund
grants \hbox{P17563-N13}, and S9607 (the second is part
of the Austrian National Research Network
"Analytic Combinatorics and Probabilistic Number Theory").}

\title[Elliptic enumeration of lattice paths]
{Elliptic enumeration of\\ nonintersecting lattice paths}

\subjclass[2000]{Primary 05A15;
Secondary 05A17, 05A19, 05E10, 11B65, 33D15, 33E20}

\keywords{nonintersecting lattice paths, elliptic weights,
elliptic hypergeometric series, Frenkel and Turaev's ${}_{10}V_9$
summation, elliptic determinant evaluations}

\newcommand{\ta}{\theta}
\newcommand{\ty}{\infty}
\newcommand{\la}{\lambda}

\begin{document}

\begin{abstract}
We enumerate lattice paths in the planar integer lattice consisting of
positively directed unit vertical and horizontal steps with respect to
a specific elliptic weight function.
The elliptic generating function of paths from a given starting point
to a given end point evaluates to an elliptic generalization of the
binomial coefficient. Convolution gives an identity equivalent to
Frenkel and Turaev's ${}_{10}V_9$ summation.
This appears to be the first combinatorial proof of the latter,
and at the same time of some important degenerate cases including
Jackson's ${}_8\phi_7$ and Dougall's ${}_7F_6$ summation.
By considering nonintersecting lattice paths we are led to
a multivariate extension of the ${}_{10}V_9$ summation which
turns out to be a special case of an identity originally
conjectured by Warnaar, later proved by Rosengren.
We conclude with discussing some future perspectives.
\end{abstract}

\maketitle

\section{Preliminaries}

\subsection{Lattice paths in $\mathbb Z^2$}

We consider lattice paths in the planar integer lattice $\mathbb Z^2$ 
consisting of unit horizontal and vertical steps in the positive direction.
Given points $u$ and $v$ in $\mathbb Z^2$, we denote the set of all
lattice paths from $u$ to $v$ by $\mathcal P(u\to v)$.
If $\mathbf u=(u_1,\dots,u_r)$ and $\mathbf v=(v_1,\dots,v_r)$ are
$r$-tuples of points, we denote the set of all $r$-tuples $(P_1,\dots,P_r)$
of paths where $P_i$ runs from $u_i$ to $v_i$, $i=1,\dots,r$, by
$\mathcal P(\mathbf u\to \mathbf v)$. A set of paths is
{\em nonintersecting} if no two paths have a point in common.
The set of all nonintersecting paths from $\mathbf u$ to $\mathbf v$
is denoted $\mathcal P_+(\mathbf u\to\mathbf v)$.
Let $w$ be a function which assigns to each horizontal edge $e$ in
$\mathbb Z^2$ a {\em weight} $w(e)$. The weight $w(P)$ of a path $P$
is defined to be the product of the weights of all its horizontal steps.
The weight $w(\mathbf P)$ of an $r$-tuple $\mathbf P=(P_1,\dots,P_r)$
of paths is defined to be the product $\prod_{i=1}^rw(P_i)$ of the weights
of all the paths in the $r$-tuple.
For any weight function $w$ defined on a set $M$, we write
$$
w(\mathcal M):=\sum_{x\in\mathcal M}w(x)
$$
for the generating function of the set $M$ with respect to the weight $w$.

For $\mathbf u=(u_1,\dots,u_r)$ and a permutation $\sigma\in\mathcal S_r$
we denote $\mathbf u_\sigma=(u_{\sigma(1)},\dots,u_{\sigma(r)})$.
We say that $\mathbf u$ is {\em} compatible to $\mathbf v$ if no
families $(P_1,\dots,P_r)$ of nonintersecting paths from
$\mathbf u_\sigma$ to $\mathbf v$ exist unless $\sigma=\epsilon$,
the identity permutation.

We need the following theorem which is a special case
(sufficient for the purposes of the present exposition)
of the Lindstr\"om--Gessel--Viennot theorem of nonintersecting
lattice paths (cf.\ \cite{L} and \cite{GV}).
\begin{theorem}\label{thnlp}
Let $\mathbf u,\mathbf v\in(\mathbb Z^2)^r$. If $\mathbf u$ is
compatible to $\mathbf v$, then 
\begin{equation}
w(\mathcal P_+(\mathbf u\to\mathbf v))=
\det_{1\le i,j\le r}w(\mathcal P(u_j\to v_i)).
\end{equation}
\end{theorem}

\subsection{Elliptic hypergeometric series}

For the following material, we refer to Chapter~11 of
Gasper and Rahman's texts \cite{GR}. Define a modified
Jacobi theta function with argument $x$ and nome $p$ by
\begin{equation}\label{mjtf}
\ta(x)=\ta (x; p):= (x; p)_\ty (p/x; p)_\ty\,,\quad\quad
\ta (x_1, \ldots, x_m) = \prod^m_{k=1} \ta (x_k),
\end{equation}
where $ x, x_1, \ldots, x_m \ne 0,\ |p| < 1,$ and $(x; p)_\ty=
\prod^\infty_{k=0}(1-x p^k)$.
We note the following useful properties of theta functions:
\begin{equation}\label{tif}
\ta(x)=-x\,\ta(1/x),
\end{equation}
\begin{equation}\label{p1id}
\ta(px)=-\frac 1x\,\ta(x),
\end{equation}
and Riemann's {\em addition formula}
\begin{equation}\label{addf}
\ta(xy,x/y,uv,u/v)-\ta(xv,x/v,uy,u/y)
=\frac uy\,\ta(yv,y/v,xu,x/u)
\end{equation}
(cf.\ \cite[p.~451, Example 5]{WW}).

Further, define a {\em theta shifted factorial} analogue of the
$q$-shifted factorial by
\begin{equation}\label{defepoch}
(a;q,p)_n = \begin{cases}
\prod^{n-1}_{k=0} \ta (aq^k),& n = 1, 2, \ldots\,,\cr
1,& n = 0,\cr
1/\prod^{-n-1}_{k=0} \ta (aq^{n+k}),& n = -1, -2, \ldots,
\end{cases}
\end{equation}
and let
\begin{equation*}
(a_1, a_2, \ldots, a_m;q, p)_n = \prod^m_{k=1} (a_k;q,p)_n,
\end{equation*}
where $a, a_1,\ldots,a_m \neq 0$.
Notice that $\ta (x;0) = 1-x$ and, hence, $(a;q, 0)_n = (a;q)_n
\big(=(a;q)_\infty/(aq^n;q)_\infty\big)$
is a $q$-{\it shifted factorial} in base $q$.
The parameters $q$ and $p$ in $(a;q,p)_n$ are called the
{\it base} and {\it nome}, respectively, and
$(a;q,p)_n$ is called the $q,p$-{\it shifted factorial}.
Observe that
\begin{equation}\label{pid}
(pa;q,p)_n=(-1)^na^{-n}q^{-\binom n2}\,(a;q,p)_n,
\end{equation}
which follows from \eqref{p1id}. 
A list of other useful identities for manipulating the
$q,p$-shifted factorials is given in \cite[Sec.~11.2]{GR}.

We call a series $\sum c_n$ an {\it elliptic hypergeometric series} if
$g(n) = c_{n+1}/c_n$ is an elliptic function of $n$ with $n$
considered as a complex variable; i.e., the function $g(x)$
is a doubly periodic meromorphic function of the complex variable $x$.
Without loss of generality, by the theory of theta functions,
we may assume that
\begin{equation*}
g(x)=\frac{\ta(a_1q^x,a_2q^x,\dots,a_{s+1}q^x;p)}
{\ta(q^{1+x},b_1q^x,\dots,b_sq^x;p)}\,z,
\end{equation*}
where the {\em elliptic balancing condition}, namely
$$a_1a_2\cdots a_{s+1}=qb_1b_2\cdots b_s,$$
holds.
If we write $q=e^{2\pi i\sigma}$, $p=e^{2\pi i\tau}$,
with complex $\sigma$, $\tau$, then $g(x)$ is indeed periodic in $x$
with periods $\sigma^{-1}$ and $\tau\sigma^{-1}$.

The general form of an elliptic hypergeometric series is thus
\begin{equation*}
{}_{s+1}E_s\!\left[\begin{matrix}a_1,\dots,a_{s+1}\\
b_1,\dots,b_s\end{matrix};q,p;z\right]
:=\sum_{k=0}^{\infty}\frac
{(a_1,a_2,\dots,a_{s+1};q,p)_k}{(q,b_1\dots,b_s;q,p)_k}z^k,
\end{equation*}
provided $a_1a_2\cdots a_{s+1}=qb_1b_2\cdots b_s$.
Here $a_1,\dots,a_r$ are the upper parameters,
$b_1,\dots,b_s$ the lower parameters, $q$ is the base,
$p$ the nome, and $z$ is the argument of the series.
For convergence reasons, one usually requires $a_{s+1}=q^{-n}$
($n$ being a nonnegative integer),
so that the sum is in fact finite.

{\em Very-well-poised elliptic hypergeometric series} are defined as
\begin{multline}\label{vwpehs}
{}_{s+1}V_s(a_1;a_6,\dots,a_{s+1};q,p;z)\\:=
{}_{s+1}E_s\!\left[\begin{matrix}{\scriptstyle a_1},\,
{\scriptstyle qa_1^{\frac 12}},{\scriptstyle -qa_1^{\frac 12}},
{\scriptstyle qa_1^{\frac 12}/p^{\frac 12}},
{\scriptstyle -qa_1^{\frac 12}p^{\frac 12}},{\scriptstyle a_6},
\dots,{\scriptstyle a_{s+1}}\\
{\scriptstyle a_1^{\frac 12}},{\scriptstyle -a_1^{\frac 12}},
{\scriptstyle a_1^{\frac 12}p^{\frac 12}},
{\scriptstyle -a_1^{\frac 12}/p^{\frac 12}},
{\scriptstyle a_1q/a_6},\dots,
{\scriptstyle a_1q/a_{s+1}}\end{matrix};q,p;-z\right]\\
=\sum_{k=0}^{\infty}\frac{\ta(a_1q^{2k};p)}{\ta(a_1;p)}
\frac{(a_1,a_6,\dots,a_{s+1};q,p)_k}
{(q,a_1q/a_6,\dots,a_1q/a_{s+1};q,p)_k}(qz)^k,
\end{multline}
where
\begin{equation*}
q^2a_6^2a_7^2\cdots a_{s+1}^2=(a_1q)^{s-5}.
\end{equation*}
It is convenient to abbreviate
\begin{equation*}
{}_{s+1}V_s(a_1;a_6,\dots,a_{s+1};q,p)
:={}_{s+1}V_s(a_1;a_6,\dots,a_{s+1};q,p;1).
\end{equation*}
Note that in \eqref{vwpehs} we have used
\begin{equation*}
\frac{\ta(aq^{2k};p)}{\ta(a;p)}=
\frac{(qa^{\frac 12},-qa^{\frac 12},qa^{\frac 12}/p^{\frac 12},
-qa^{\frac 12}p^{\frac 12};q,p)_k}
{(a^{\frac 12},-a^{\frac 12},
a^{\frac 12}p^{\frac 12},-a^{\frac 12}/p^{\frac 12};q,p)_k}
(-q)^{-k},
\end{equation*}
which shows that in the elliptic case the number of
pairs of numerator and denominator paramters
involved in the construction of the {\em very-well-poised term} is {\em four}
(whereas in the basic case this number is {\em two},
in the ordinary case only {\em one}).

The above definitions for ${}_{s+1}E_{s}$ and ${}_{s+1}V_{s}$ series
are due to Spiridonov~\cite{Sp}, see \cite[Ch.~11]{GR}.

In their study of elliptic $6j$ symbols (which are elliptic solutions
of the Yang--Baxter equation found by Baxter~\cite{B} and Date et
al.~\cite{DJKMO}), Frenkel and Turaev~\cite{FT}
came across the following ${}_{12}V_{11}$ transformation:
\begin{multline}\label{12V11}
{}_{12}V_{11}(a;b,c,d,e,f,\lambda a q^{n+1}/ef,q^{-n};q,p)
=\frac {(aq,aq/ef,\lambda q/e,\lambda q/f;q,p)_n}
{(aq/e,aq/f,\lambda q/ef,\lambda q;q,p)_n}\\\times
{}_{12}V_{11}(\lambda;\lambda b/a,\lambda c/a,\lambda d/a,e,f,
\lambda a q^{n+1}/ef,q^{-n};q,p),
\end{multline}
where $\lambda=a^2q/bcd$. This is an extension of Bailey's
very-well-poised ${}_{10}\phi_9$ transformation~\cite[Eq.~(2.9.1)]{GR},
to which it reduces when $p=0$.

The ${}_{12}V_{11}$ transformation in \eqref{12V11} appeared as 
a consequence of the tetrahedral symmetry of the elliptic $6j$ symbols.
Frenkel and Turaev's transformation contains as a special case
the following summation formula,
\begin{equation}
{}_{10}V_9(a;b,c,d,e,q^{-n};q,p)
=\frac {(aq,aq/bc,aq/bd,aq/cd;q,p)_n}
{(aq/b,aq/c,aq/d,aq/bcd;q,p)_n},
\end{equation}
where $a^2q^{n+1}=bcde$, see also \eqref{10V9sum}.
The $_{10}V_9$ summation is an elliptic analogue of Jackson's
$_8\phi_7$ summation formula~\cite[Eq.~(2.6.2)]{GR}
(or of Dougall's $_7F_6$ summation formula~\cite[Eq.~(2.1.6)]{GR}).
A striking feature of elliptic hypergeometric series is that
already the simplest identities involve many parameters.
The fundamental identity at the ``bottom'' of the hierarchy of
identities for elliptic hypergeometric series is the $_{10}V_9$ summation.
When keeping the nome $p$ arbitrary (while $|p|<1$) there is no way
to specialize (for the sake of obtaining lower order identities)
any of the free parameters of an elliptic hypergeometric series
in form of a limit tending to zero or infinity, due to the issue
of convergence. For the same reason, elliptic hypergeometric series
are only well-defined as complex functions if they are terminating
(i.e., the sums are finite).
See Gasper and Rahman's texts \cite[Ch.~11]{GR} for more details.

The outline of the remaining sections of this paper is as follows:
In Section~\ref{secellen} we introduce a specific elliptic weight
function, composed of appropriately chosen products of theta functions.
Using this weight, we then compute the elliptic generating function
of paths from a given starting point to a given end point.
The result simplifies, by virtue of Riemann's addition formula for
theta functions and induction, to closed form, namely
to an elliptic generalization of the binomial coefficient.
By convolution we readily obtain an identity equivalent to
Frenkel and Turaev's ${}_{10}V_9$ summation.
This appears to be the first combinatorial proof of this important
summation (fundamental to the theory of elliptic hypergeometric series),
and at the same time of some important degenerate
cases including Jackson's ${}_8\phi_7$ and Dougall's ${}_7F_6$
summation, both fundamental to the respective theories of basic
and ordinary hypergeometric series.
We then turn to nonintersecting lattice paths in $\mathbb Z^2$ where,
using the Lindstr\"om--Gessel--Viennot theorem combined with an
elliptic determinant evaluation by Warnaar, we compute
the elliptic generating function of selected families of paths
with given starting points and end points.
Here convolution gives a multivariate extension of the
${}_{10}V_9$ summation, see Section~\ref{secms}, which turns out
to be a special case of an identity originally conjectured by
Warnaar, later proved by Rosengren. We also display a more
general multivariate ${}_{12}V_{11}$ transformation
(being a special case of an identity originally conjectured by Warnaar,
later proved by Rains, and, independently, by Coskun and Gustafson),
which we strongly believe can be established by the methods of
this paper, which however we were so far unable to accomplish.
We conclude in Section~\ref{secfp} with discussing some future
perspectives, in particular, concerning the elliptic enumeration
of tableaux and plane partitions, a variant of elliptic Schur
functions, other weight functions, and the commencement of general
research in ``elliptic combinatorics''.

\section{Elliptic enumeration of lattice paths}\label{secellen}

The identity responsible for $q$-calculus to ``work'' is
the simple factorization
\begin{equation}\label{sid}
q^k-q^{k+1}=(1-q)q^k.
\end{equation}
This (almost embarrassingly simple) identity underlies not only
$q$-integration (cf.\ \cite[Eq.~(2.12)]{A}),
but also the recursion(s) for the $q$-binomial coefficient
(see \eqref{recqbin} at the end of this section).
As $q$-binomial coefficients can be combinatorially interpreted
as generating functions of lattice paths in $\mathbb Z^2$
(from a given starting point to a given end point),
one may wonder whether any suitable generalization of \eqref{sid}
would give rise to a corresponding extension of
$q$-binomial coefficients with meaningful combinatorial
interpretation.
Indeed, by using the much more general identity \eqref{addf},
rather than \eqref{sid}, as the underlying three term relation,
we obtain such an extension. In particular,
we shall be considering {\em elliptic binomial coefficients},
resulting from the enumeration of lattice paths with respect to
{\em elliptic weights}.
The expressions and series occurring in our study belong to the world
of elliptic hypergeometric series, which we just introduced
in the previous section. 

The most important ingredient for this analysis to work out
is the particular ``clever'' choice of weight function in \eqref{ellw}.
This choice was made, on one hand, by matching the general indefinite sum
\eqref{sindefsum} with the known indefinite sum in \eqref{windefsum},
such that induction can be applied
(with appeal to the three term relation \eqref{addf},
actually a special case of \eqref{windefsum}). One the other hand,
factorization of the elliptic binomial coefficient
$w(\mathcal P((l,k)\to(n,m)))$ was sought in general, in particular
also when $(l,k)\neq (0,0)$. Once the right choice of
weight function is made, everything becomes easy and a matter
of pure verification. Nevertheless, at the conceptual level things
remain interesting (and non-trivial).
For instance, the elliptic binomial coefficient
$w(\mathcal P((l,k)\to(n,m)))$ indeed depends on $l,k,n,m$
(besides other parameters), and is not a mere multiple of
$w(\mathcal P((0,0)\to(n-l,m-k)))$, contrary to the basic (``$q$'')
or classical case. 

Let $a,b,q,p$ be arbitrary (complex) parameters with $a,b,q\neq 0$
and $|p|<1$. We define the (``standard'') {\em elliptic weight function}
on horizontal edges $(n-1,m)\to(n,m)$ of $\mathbb Z^2$ as follows.
\begin{multline}\label{ellw}
w(n,m)=w(n,m;a,b;q,p)\\:=
\frac{\ta(aq^{n+2m},bq^{2n},bq^{2n-1},aq^{1-n}/b,aq^{-n}/b)}
{\ta(aq^n,bq^{2n+m},bq^{2n+m-1},aq^{1+m-n}/b,aq^{m-n}/b)}q^m.
\end{multline}

Our terminology is perfectly justified as the weight function
defined in \eqref{ellw} is indeed {\em elliptic} (i.e., doubly periodic
meromorphic), even independently in each of $\log_qa$, $\log_qb$,
$n$ and $m$ (viewed as complex parameters). 
If we write $q=e^{2\pi i\sigma}$, $p=e^{2\pi i\tau}$, $a=q^\alpha$ and
$b=q^\beta$ with complex $\sigma$, $\tau$, $\alpha$ and
$\beta$, then the weight $w(n,m)$ is clearly periodic in $\alpha$
with period $\sigma^{-1}$. A simple calculation involving \eqref{pid}
further shows that $w(n,m)$ is also periodic in $\alpha$
with period $\tau\sigma^{-1}$ (the latter means that $w(n,m)$
is invariant with respect to $a\mapsto pa$).
The same applies to $w(n,m)$ viewed as a function in $\beta$
(or $n$ or $m$) with the same two periods $\sigma^{-1}$ and $\tau\sigma^{-1}$.
Spiridonov~\cite{Sp} calls expressions such as \eqref{ellw} where
all free parameters have equal periods of double periodicity
{\em totally elliptic}. In this respect we can also refer to
\eqref{ellw} as a totally elliptic weight.

For $p=0$ \eqref{ellw} reduces to
\begin{multline}\label{vwpw}
w(n,m;a,b;q,0)\\=
\frac{(1-aq^{n+2m})(1-bq^{2n})(1-bq^{2n-1})(1-aq^{1-n}/b)(1-aq^{-n}/b)}
{(1-aq^n)(1-bq^{2n+m})(1-bq^{2n+m-1})(1-aq^{1+m-n}/b)(1-aq^{m-n}/b)}q^m.
\end{multline}
If we further let $a\to 0$ and then $b\to 0$ (in this order; or take $b\to 0$
and then $a\to\infty$) this reduces to the standard $q$-weight $q^m$
(counting the height of, or the area below, the horizontal edge
$(n-1,m)\to(n,m)$).

By an {\em elliptic generating function} we mean, of course, a generating
function with respect to an elliptic weight function (and in
particular, we shall always take the weight defined in \eqref{ellw}
unless stated otherwise). It is clear that an elliptic generating
function is elliptic as a function in its free parameters.

The particular choice of our elliptic weight in \eqref{ellw} is
justified by the following very nice result.

\begin{theorem}\label{mainres}
Let $l,k,n,m$ be four integers with $n-l+m-k\ge 0$.
The elliptic generating function of paths running from
$(l,k)$ to $(n,m)$ is
\begin{multline}\label{ellbc}
w(\mathcal P((l,k)\to(n,m)))=
\frac{(q^{1+n-l},aq^{1+n+2k},bq^{1+n+k+l},aq^{1+k-n}/b;q,p)_{m-k}}
{(q,aq^{1+l+2k},bq^{1+2n+k},aq^{1+k-l}/b;q,p)_{m-k}}\\\times
\frac{(aq^{1+l+2k},aq^{1-n}/b,aq^{-n}/b;q,p)_{n-l}}
{(aq^{1+l},aq^{1+k-n}/b,aq^{k-n}/b;q,p)_{n-l}}
\frac{(bq^{1+2l};q,p)_{2n-2l}}{(bq^{1+k+2l};q,p)_{2n-2l}}q^{(n-l)k}.
\end{multline}
\end{theorem}

\begin{proof}
First, if $k>m$ (there is no path in this case), the expression in
\eqref{ellbc} vanishes due to the factor $(q;q,p)_{m-k}^{-1}$.
On the other hand, if $m\ge k$ but $l>n$ (again there is no path)
the expression vanishes due to the factor $(q^{1+n-l};q,p)_{m-k}$
since $n-l+m-k\ge 0$. We may therefore assume, besides
$n-l+m-k\ge 0$, that $n\ge l$ and $m\ge k$. The statement is now
readily proved by induction on $n-l+m-k$.
For $n=l$ one has $w(\mathcal P((l,k)\to(l,m)))=1$ as desired.
For $m=k$ one readily verifies $w(\mathcal P((l,k)\to(n,k)))=
\prod_{i=l+1}^nw(i,k)$. (In both cases there is just one path.)
Next assume $n>l$ and $m>k$. We are done if we can verify
the recursion
\begin{multline}\label{aid}
w(\mathcal P((l,k)\to(n,m)))\\=
w(\mathcal P((l,k)\to(n,m-1)))+w(\mathcal P((l,k)\to(n-1,m)))\,w(n,m).
\end{multline}
(The final step of a path is either vertical or horizontal.)
However, after cancellation of common factors this reduces
to the addition formula \eqref{addf}.
\end{proof}

Aside from the recursion \eqref{aid},
we also (automatically) have
\begin{multline}\label{aid2}
w(\mathcal P((l,k)\to(n,m)))\\=
w(\mathcal P((l,k+1)\to(n,m)))+w(l+1,k)\,w(\mathcal P((l+1,k)\to(n,m))).
\end{multline}
(The first step of a path is either vertical or horizontal.)
In the limit $p\to 0$, $a\to 0$, $b\to 0$ (in this order), the
recursions \eqref{aid} and \eqref{aid2} reduce to
\begin{multline*}
\left[\begin{matrix}n-l+m-k\\n-l\end{matrix}\right]_qq^{(n-l)k}\\
=\left[\begin{matrix}n-l+m-k-1\\n-l\end{matrix}\right]_qq^{(n-l)k}
+\left[\begin{matrix}n-l+m-k-1\\n-l-1\end{matrix}\right]_qq^{(n-l-1)k+m}
\end{multline*}
and
\begin{multline*}
\left[\begin{matrix}n-l+m-k\\n-l\end{matrix}\right]_qq^{(n-l)k}\\
=\left[\begin{matrix}n-l+m-k-1\\n-l\end{matrix}\right]_qq^{(n-l)(k+1)}
+\left[\begin{matrix}n-l+m-k-1\\n-l-1\end{matrix}\right]_qq^{(n-l-1)k+k},
\end{multline*}
respectively,
where
\begin{equation}\label{qbc}
\left[\begin{matrix}n\\k\end{matrix}\right]_q:=
\frac{(q;q)_n}{(q;q)_k(q;q)_{n-k}}
\end{equation}
is the {\em $q$-binomial coefficient}, defined for nonnegative
integers $n,k$ with $n\ge k$.
This pair of recursions is of course equivalent to the well-known pair
\begin{equation}\label{recqbin}
\left[\begin{matrix}n\\k\end{matrix}\right]_q
=\left[\begin{matrix}n-1\\k\end{matrix}\right]_q
+\left[\begin{matrix}n-1\\k-1\end{matrix}\right]_qq^{n-k},\qquad
\left[\begin{matrix}n\\k\end{matrix}\right]_q
=\left[\begin{matrix}n-1\\k\end{matrix}\right]_qq^k
+\left[\begin{matrix}n-1\\k-1\end{matrix}\right]_q.
\end{equation}

We may therefore refer to the factored expression in \eqref{ellbc}
as an {\em elliptic binomial coefficient} (which should not be
confused with the much simpler definition given in
\cite[Eq.~(11.2.61)]{GR} which is a straightforward theta shifted
factorial extension of \eqref{qbc} but actually {\em not} elliptic).
In fact, it is not difficult to see that the expression in
\eqref{ellbc} is totally elliptic, i.e.\
elliptic in each of $\log_qa$, $\log_qb$, $l$, $k$, $n$ and $m$
(viewed as complex parameters) which again fully justifies
the notion ``elliptic''.

\begin{remark}\label{remsh}
Consider the two parameter extension of \eqref{ellw} defined by
\begin{equation}\label{ellw2}
w_{(s,t)}(n,m;a,b;q,p):=w(n,m;aq^{s+2t},bq^{2s+t};q,p).
\end{equation}
Clearly, $w_{(0,0)}(n,m)=w(n,m)$.
A simple calculation reveals that
\begin{equation}
w(n+s,m+t)=w_{(s,0)}(n,t)\,w_{(s,t)}(n,m).
\end{equation}
This notation is useful for dealing with shifted paths.
In terms of generating functions we have
\begin{multline}\label{bid}
w(\mathcal P((l+s,k+t)\to(n+s,m+t)))\\=
w_{(s,0)}(\mathcal P((l,t)\to(n,t)))\,
w_{(s,t)}(\mathcal P((l,k)\to(n,m))),
\end{multline}
which is readily verified using Theorem~\ref{mainres}.

Other useful properties of the weight function in \eqref{ellw} are
\begin{align}
w(n,m;a,b;q,p)&=w(-n,-m;a^{-1},qb^{-1};q,p)\\
&=w(n,m;a^{-1},b^{-1};q^{-1},p).
\end{align}
Furthermore, invoking Theorem~\ref{mainres} one easily verifies
\begin{multline}\label{cid}
w(\mathcal P((l,k)\to(n,m));a,b;q,p)\\=
w(\mathcal P((-1-n,-m)\to(-1-l,-k));a^{-1},qb^{-1};q,p).
\end{multline}
\end{remark}

\subsection{Immediate consequences}

Let us consider the elliptic generating function of lattice paths
in $\mathbb Z^2$ from $(0,0)$ to $(n,m)$. (In what follows, there
is in fact no loss of generality in choosing the starting point
to be the origin.) We may distinguish the paths according to the
height of the last step.
This gives the simple identity
\begin{equation}\label{sindefsum}
w(\mathcal P((0,0)\to(n,m)))=
\sum_{k=0}^mw(\mathcal P((0,0)\to(n-1,k)))\,w(n,k).
\end{equation}
In explicit terms, this is
\begin{multline*}
\frac{(q^{1+n},aq^{1+n},bq^{1+n},aq^{1-n}/b;q,p)_{m}}
{(q,aq,bq^{1+2n},aq/b;q,p)_{m}}=\\
\sum_{k=0}^m\frac{(q^{n},aq^{n},bq^{n},aq^{2-n}/b;q,p)_{k}\,
\ta(aq^{n+2k},bq^{2n},bq^{2n-1},aq^{1-n}/b,aq^{-n}/b)}
{(q,aq,bq^{2n-1},aq/b;q,p)_{k}\,
\ta(aq^n,bq^{2n+k},bq^{2n+k-1},aq^{1+k-n}/b,aq^{k-n}/b)}q^k,
\end{multline*}
which, after simplifying the summand, is
\begin{equation}\label{indefsum}
\frac{(q^{1+n},aq^{1+n},bq^{1+n},aq^{1-n}/b;q,p)_{m}}
{(q,aq,bq^{1+2n},aq/b;q,p)_{m}}
=\sum_{k=0}^m\frac{\ta(aq^{n+2k})(aq^{n},q^{n},bq^{n},aq^{-n}/b;q,p)_{k}}
{\ta(aq^{n})(q,aq,aq/b,bq^{1+2n};q,p)_{k}}q^k.
\end{equation}
By analytic continuation to replace $q^n$ by an arbitrary complex
parameter (\eqref{indefsum} is true for all $n\ge 0$, etc.;
see Warnaar~\cite[Proof of Thms.\ 4.7--4.9]{Wa} for a typical
application of the identity theorem in the elliptic setting)
and substitution of variables, one gets the indefinite summation
\begin{equation}\label{windefsum}
\frac{(aq,bq,cq,aq/bc;q,p)_{m}}
{(q,aq/b,aq/c,bcq;q,p)_{m}}
=\sum_{k=0}^m\frac{\ta(aq^{2k})(a,b,c,a/bc;q,p)_{k}}
{\ta(a)(q,aq/b,aq/c,bcq;q,p)_{k}}q^k
\end{equation}
(cf.\ \cite[Eq.~(11.4.10)]{GR}).

More generally, for a fixed $l$, $1\le l\le n$, we may distinguish
paths running from $(0,0)$ to $(n,m)$ by the height $k$ they have
when they first reach a point on the vertical line $x=l$ (right after
the horizontal step $(l-1,k)\to(l,k)$). This refined enumeration reads,
in terms of elliptic generating functions,
\begin{multline}\label{ft1}
w(\mathcal P((0,0)\to(n,m)))\\=
\sum_{k=0}^mw(\mathcal P((0,0)\to(l-1,k)))\,w(l,k)\,
w(\mathcal P((l,k)\to(n,m))).
\end{multline}
Explicitly, this is (after some simplifictions)
\begin{multline}
\frac{(q^{1+n},aq^{1+l},bq^{1+n},aq^{1-l}/b;q,p)_{m}}
{(q^{1+n-l},aq,bq^{1+n+l},aq/b;q,p)_{m}}\\
=\sum_{k=0}^m\frac{\ta(aq^{l+2k})
(aq^{l},bq^{l},q^{l},aq^{-n}/b,aq^{1+n+m},q^{-m};q,p)_{k}}
{\ta(aq^{l})(q,aq/b,aq,bq^{1+n+l},q^{l-n-m},aq^{1+l+m};q,p)_{k}}q^k,
\end{multline}
which after analytic continuation (first to replace $q^n$, then $q^l$,
by complex parameters) and substitution of variables becomes
\begin{multline}\label{10V9sum}
\frac{(aq,aq/bc,aq/bd,aq/cd;q,p)_{m}}
{(aq/b,aq/c,aq/d,aq/bcd;q,p)_{m}}\\
=\sum_{k=0}^m\frac{\ta(aq^{2k})
(a,b,c,d,a^2q^{1+m}/bcd,q^{-m};q,p)_{k}}
{\ta(a)(q,aq/b,aq/c,aq/d,bcdq^{-m}/a,aq^{1+m};q,p)_{k}}q^k,
\end{multline}
The result is {\em Frenkel and Turaev's ${}_{10}V_9$ summation}
(\cite{FT}; cf.\ \cite[Eq.~(11.4.1)]{GR}), the elliptic extension of
Jackson's very-well-poised balanced ${}_8\phi_7$ summation
(cf.\ \cite[Eq.~(2.6.2)]{GR}), the latter of which is a
$q$-analogue of Dougall's ${}_7F_6$ summation theorem.
Of course, the $p\to 0$ limit case of the above analysis
(using the weight function in \eqref{vwpw}) reduces to a
proof of Jackson's ${}_8\phi_7$ summation.
On the other hand, the $p\to 0$, $a\to 0$ limit case of this analysis,
with the weight function
\begin{equation}
w(n,m;0,b;q,0)\\=
\frac{(1-bq^{2n})(1-bq^{2n-1})}
{(1-bq^{2n+m})(1-bq^{2n+m-1})}q^m,
\end{equation}
yields the $q$-Pfaff--Saalsch\"utz summation for a balanced
terminating $_3\phi_2$ series (cf.\ \cite[Eq.~(1.7.2)]{GR}).
(A completely different combinatorial proof of the $_3\phi_2$ summation
was given by Zeilberger~\cite{Z}.) If one further lets (in addition to
$p\to 0$ and $a\to 0$) $b\to 0$, where one considers the standard $q$-weight,
the above analysis yields, as is well-known, the $q$-Chu--Vandermonde
summation (cf.\ \cite[Eq.~(1.5.3)]{GR}).

We briefly sketch two other ways how to obtain the ${}_{10}V_9$ sum
from Theorem~\ref{mainres} by convolution (and analytic continuation).
For a fixed $k$, $1\le k\le m$, we may distinguish
paths running from $(0,0)$ to $(n,m)$ by the abscissa $l$ they have
when they first reach a point on the horizontal line $y=k$ (right after
the vertical step $(l,k-1)\to(l,k)$). This refined enumeration reads,
in terms of elliptic generating functions,
\begin{equation}\label{ft2}
w(\mathcal P((0,0)\to(n,m)))=
\sum_{l=0}^mw(\mathcal P((0,0)\to(l,k-1)))\,w(\mathcal P((l,k)\to(n,m))).
\end{equation}
On the other hand, we may also fix an antidiagonal running through
$(k,0)$ and $(0,k)$, $0<k<n+m$. We can then distinguish
paths running from $(0,0)$ to $(n,m)$ by where they cut
the antidiagonal. This refined enumeration reads,
in terms of elliptic generating functions,
\begin{equation}\label{ft3}
w(\mathcal P((0,0)\to(n,m)))=
\sum_{l=0}^{\min(k,n)}w(\mathcal P((0,0)\to(l,k-l)))\,
w(\mathcal P((l,k-l)\to(n,m))).
\end{equation}
The last two identities both constitute, when written out
explicitly using Theorem~\ref{mainres}, variants of Frenkel and
Turaev's ${}_{10}V_9$ summation (like \eqref{ft1}) both
of which can be extended to
\eqref{10V9sum} by analytic continuation.

\subsection{Determinant evaluations and elliptic generating functions
for nonintersecting lattice paths}

For obtaining explicit results the following determinant evaluation,
taken from \cite[Cor.~5.4]{Wa}, is crucial.
\begin{lemma}[Warnaar]\label{lemdet1}
Let $A$, $B$, $C$, and $X_1,\dots,X_r$ be indeterminate. Then there holds
\begin{multline}\label{lemdet1gl}
\det_{1\le i,j\le r}\left(
\frac{(AX_i,AC/X_i;q,p)_{r-j}}
{(BX_i,BC/X_i;q,p)_{r-j}}\right)
=A^{\binom r2}q^{\binom r3}
\prod_{1\le i<j\le r}X_j\,\ta(X_i/X_j,C/X_iX_j)\\\times
\prod_{i=1}^{r}\frac {(B/A,ABCq^{2r-2i};q,p)_{i-1}}
{(BX_i,BC/X_i;q,p)_{r-1}}.
\end{multline}
\end{lemma}

As a consequence of Theorem~\ref{thnlp} and Lemma~\ref{lemdet1},
we have the following explicit formulae which generalize
Theorem~\ref{mainres}:
\begin{proposition}\label{prop}
(a) Let $l,k,n$, $m_1,\dots,m_r$ be integers such that
$m_1\ge m_2\ge\dots\ge m_r$ and $n-l+m_i-k\ge 0$ for all $i=1,\dots,r$.
Then the elliptic generating function for nonintersecting
lattice paths with starting points $(l+i,k-i)$ and end points
$(n,m_i)$, $i=1,\dots,r$, is
\begin{multline}\label{det1}
\det_{1\le i,j,\le r}\big(w(\mathcal P((l+j,k-j)\to(n,m_i)))\big)\\
=q^{3\binom{r+1}3+\binom{r+2}3+r(n-l)k-(n-l)\binom{r+1}2-r^2k+
\sum_{i=1}^r(i-1)m_i}\\\times
\prod_{1\le i<j\le r}\ta(q^{m_i-m_j},aq^{1+n+m_i+m_j})\\\times
\prod_{i=1}^r\frac{(q^{1+n-l-i};q,p)_{m_i-k+i}
(aq^{1+n+2k-r-i};q,p)_{m_i-k+i}(aq^{1+l+2k-i};q,p)_{n-l-r}}
{(q;q,p)_{m_i-k+r}
(aq^{1+l+2k-i};q,p)_{m_i-k+i}(aq^{1+l+i};q,p)_{n-l-i}}\\\times
\prod_{i=1}^r\frac{(bq^{2+n+k+l-i};q,p)_{m_i-k+i}(bq^{1+2l+2i};q,p)_{2n-2l-2i}}
{(bq^{1+2n+k-i};q,p)_{m_i-k+i}
(bq^{1+2l+k+i};q,p)_{2n-2l-2i}}\\\times
\prod_{i=1}^r\frac{(aq^{1+k-n-i}/b;q,p)_{m_i-k+i}
(aq^{1-n}/b,aq^{-n}/b;q,p)_{n-l-i}}
{(aq^{k-l-i}/b;q,p)_{m_i-k+i}(aq^{1+k-n-i}/b,aq^{k-n-i}/b;q,p)_{n-l-i}}.
\end{multline}

(b) Let $l,k,m$, $n_1,\dots,n_r$, be integers such that
$n_1\le n_2\le\dots\le n_r$ and $n_i-l+m-k\ge 0$ for all $i=1,\dots,r$.
Then the elliptic generating function for nonintersecting
lattice paths with starting points $(l+i,k-i)$ and end points
$(n_i,m)$, $i=1,\dots,r$, is
\begin{multline}\label{det2}
\det_{1\le i,j,\le r}\big(w(\mathcal P((l+j,k-j)\to(n_i,m)))\big)\\
=q^{2\binom{r+1}3+(l-k+1)\binom{r+1}2-rlk+\sum_{i=1}^r(k-i)n_i}
\prod_{1\le i<j\le r}
\ta(q^{n_j-n_i},bq^{1+m+n_i+n_j})\\\times
\prod_{i=1}^r\frac{(q^{n_i-l};q,p)_{m-k+1}
(aq^{n_i+2k-i};q,p)_{m-k+1-r+i}
(aq^{l+2k};q,p)_{n_i-l-i}}
{(q;q,p)_{m-k+i}(aq^{l+2k};q,p)_{m-k+1-r+i}(aq^{1+l+i};q,p)_{n_i-l-i}}\\\times
\prod_{i=1}^r\frac{(bq^{1+n_i+k+l};q,p)_{m-k+1}(bq^{1+2l+2i};q,p)_{2n_i-2l-2i}}
{(bq^{1+2n_i+k-i};q,p)_{m-k+i}
(bq^{1+2l+k+i};q,p)_{2n_i-2l-2i}}\\\times
\prod_{i=1}^r\frac{(aq^{k-n_i}/b;q,p)_{m-k+1}
(aq^{1-n_i}/b,aq^{-n_i}/b;q,p)_{n_i-l-i}}
{(aq^{k-l-i}/b;q,p)_{m-k+1}(aq^{k-n_i}/b;q,p)_{n_i-l+1-2i}
(aq^{k-r-n_i}/b;q,p)_{n_i-l+r-2i}}.
\end{multline}

(c) Let $l,k,m$, $n_1,\dots,n_r$ be integers such that
$n_1\le n_2\le\dots\le n_r$ and $m-l-k\ge 0$.
Then the elliptic generating function for nonintersecting
lattice paths with starting points $(l+i,k-i)$ and end points
$(n_i,m-n_i)$, $i=1,\dots,r$, is
\begin{multline}\label{det3}
\det_{1\le i,j,\le r}\big(w(\mathcal P((l+j,k-j)\to(n_i,m-n_i)))\big)\\
=q^{2\binom{r+1}3+(l-k+1)\binom{r+1}2-rlk+
\sum_{i=1}^r(k-i)n_i}\prod_{1\le i<j\le r}
\ta(q^{n_j-n_i},aq^{m-n_i-n_j}/b)\\\times
\prod_{i=1}^r\frac{(q^{n_i-l};q,p)_{m-n_i-k+i}
(aq^{n_i+2k-i};q,p)_{m-n_i-k+1}(aq^{l+2k};q,p)_{n_i-l-i}}
{(q;q,p)_{m-n_i-k+r}
(aq^{l+2k};q,p)_{m-n_i-k+1}(aq^{1+l+i};q,p)_{n_i-l-i}}\\\times
\prod_{i=1}^r\frac{(bq^{1+n_i+k+l};q,p)_{m-n_i-k+i}
(bq^{1+2l+2i};q,p)_{2n_i-2l-2i}}
{(bq^{1+2n_i+k-i};q,p)_{m-n_i-k+i}
(bq^{1+2l+k+i};q,p)_{2n_i-2l-2i}}\\\times
\prod_{i=1}^r\frac{(aq^{1+k-n_i-i}/b;q,p)_{m-n_i-k+i}}
{(aq^{k-l-i}/b;q,p)_{m-n_i-k+i}}\\\times
\prod_{i=1}^r\frac{(aq^{1-n_i}/b,aq^{-n_i}/b;q,p)_{n_i-l-i}}
{(aq^{1+k-n_i-i}/b;q,p)_{n_i-l-i}
(aq^{k-r-n_i}/b;q,p)_{n_i-l+r-2i}}.
\end{multline}

(d) Let $l,n,m$, $k_1,\dots,k_r$ be integers such that
$k_1\ge k_2\ge\dots\ge k_r$ and $n-l+m-k_i\ge 0$ for all $i=1,\dots,r$.
Then the elliptic generating function for nonintersecting
lattice paths with starting points $(l,k_i)$ and end points
$(n+i,m-i)$, $i=1,\dots,r$, is
\begin{multline}\label{det4}
\det_{1\le i,j,\le r}\big(w(\mathcal P((l,k_j)\to(n+i,m-i)))\big)\\
=q^{\sum_{i=1}^r(n-l+i)k_i}\prod_{1\le i<j\le r}
\ta(q^{k_i-k_j},aq^{l+k_i+k_j})\\\times
\prod_{i=1}^r\frac{(q^{1+n+i-l};q,p)_{m-k_i-i}
(aq^{1+n+2k_i};q,p)_{m-k_i}(aq^{1+l+2k_i};q,p)_{n-l}}
{(q;q,p)_{m-k_i-1}(aq^{1+l+2k_i};q,p)_{m-k_i-1}(aq^{1+l};q,p)_{n-l+i}}\\\times
\prod_{i=1}^r\frac{(bq^{1+n+k_i+l+r};q,p)_{m-k_i-r-1+i}
(bq^{1+2l};q,p)_{2n-2l+2i}}
{(bq^{1+2n+k_i};q,p)_{m-k_i+i}
(bq^{1+2l+k_i};q,p)_{2n-2l}}\\\times
\prod_{i=1}^r\frac{(aq^{1+k_i-n}/b;q,p)_{m-k_i-i-1}
(aq^{1-n-i}/b,aq^{-n-i}/b;q,p)_{n-l+i}}
{(aq^{1+k_i-l}/b;q,p)_{m-k-i}(aq^{1+k_i-n}/b;q,p)_{n-l}
(aq^{k_i-r-n}/b;q,p)_{n-l+r}}.
\end{multline}

(e) Let $k,n,m$, $l_1,\dots,l_r$ be integers such that
$l_1\le l_2\le\dots\le l_r$ and $n-l_i+m-k\ge 0$ for all $i=1,\dots,r$.
Then the elliptic generating function for nonintersecting
lattice paths with starting points $(l_i,k)$ and end points
$(n+i,m-i)$, $i=1,\dots,r$, is
\begin{multline}\label{det5}
\det_{1\le i,j,\le r}\big(w(\mathcal P((l_j,k)\to(n+i,m-i)))\big)\\
=q^{(n+r+k)\binom r2+(n+1)rk-\sum_{i=1}^r(k+i-1)l_i}\prod_{1\le i<j\le r}
\ta(q^{l_j-l_i},bq^{k+l_i+l_j})\\\times
\prod_{i=1}^r\frac{(q^{1+n+r-l_i};q,p)_{m-k-r}
(aq^{1+n+2k+i};q,p)_{m-k-1}(aq^{1+l_i+2k};q,p)_{n+i-l_i}}
{(q;q,p)_{m-k-i}(aq^{1+l_i+2k};q,p)_{m-k-1}(aq^{1+l_i};q,p)_{n+i-l_i}}\\\times
\prod_{i=1}^r\frac{(bq^{1+n+k+r+l_i};q,p)_{m-k-r}
(bq^{1+2l_i};q,p)_{2n+2i-2l_i}}
{(bq^{1+2n+k+2_i};q,p)_{m-k-i}(bq^{1+k+2l_i};q,p)_{2n+2i-2l_i}}\\\times
\prod_{i=1}^r\frac{(aq^{1+k-n-i}/b;q,p)_{m-k-1}
(aq^{1-n-i}/b,aq^{-n-i}/b;q,p)_{n+i-l_i}}
{(aq^{1+k-l_i}/b;q,p)_{m-k-1}(aq^{1+k-n-i}/b,aq^{k-n-i}/b;q,p)_{n+i-l_i}}.
\end{multline}

(f) Let $k,n,m$, $l_1,\dots,l_r$ be integers such that
$l_1\le l_2\le\dots\le l_r$ and $n+m-k\ge 0$.
Then the elliptic generating function for nonintersecting
lattice paths with starting points $(l_i,k-l_i)$ and end points
$(n+i,m-i)$, $i=1,\dots,r$, is
\begin{multline}\label{det6}
\det_{1\le i,j,\le r}\big(w(\mathcal P((l_j,k-l_j)\to(n+i,m-i)))\big)\\
=q^{k\binom{r+1}2+rnk-\sum_{i=1}^r(n+k+i-l_i)l_i}\prod_{1\le i<j\le r}
\ta(q^{l_j-l_i},aq^{k-l_i-l_j}/b)\\\times
\prod_{i=1}^r\frac{(q^{1+n+r-l_i};q,p)_{m-k-r+l_i+i-1}
(aq^{1+n+2k-2l_i};q,p)_{m-k+l_i}(aq^{1+2k-l_i};q,p)_{n-l_i}}
{(q;q,p)_{m-k+l_i-1}(aq^{1+2k-l_i};q,p)_{m-k+l_i-i}
(aq^{1+l_i};q,p)_{n+i-l_i}}\\\times
\prod_{i=1}^r\frac{(bq^{1+n+k+i};q,p)_{m-k+l_i-i}
(bq^{1+2l_i};q,p)_{2n+2i-2l_i}}
{(bq^{1+2n+k-l_i};q,p)_{m-k+l_i+i}(bq^{1+k+l_i};q,p)_{2n-2l_i}}\\\times
\prod_{i=1}^r\frac{(aq^{1+k-n-l_i}/b;q,p)_{m-k+l_i-i-1}}
{(aq^{1+k-2l_i}/b;q,p)_{m-k+l_i-1}}\\\times
\prod_{i=1}^r\frac{(aq^{1-n-i}/b,aq^{-n-i}/b;q,p)_{n+i-l_i}}
{(aq^{1+k-n-l_i}/b;q,p)_{n-l_i}(aq^{k-n-r-l_i}/b;q,p)_{n+r-l_i}}.
\end{multline}

\end{proposition}

\begin{remark}
In Proposition~\ref{prop} we are considering generating functions
for families of nonintersecting lattice paths where the set of
starting points or end points are consecutive points on an antidiagonal
parallel to $x+y=c$, for an integer $c$, such as $(l+i,c-l-i)$.
What happens if, say, the starting points are instead considered
to be consecutive points on a {\em horizontal} (resp.\ {\em vertical})
line, such as $(l+i,k)$ (resp.\ $(l,k-i)$), $i=1,\dots,r$?
The answer is that the computation of the generating function is
then readily reduced to the previous case where the starting points
are consecutive points on an antidiagonal, namely $(l+i,k+r-i)$
(resp.\ $(l+i-1,k-i)$), $i=1,\dots,r$. (We thank Christian Krattenthaler
for reminding us of this simple fact; during the preparations of this
paper, we had namely computed these other determinants separately
and were originally planning to include them explicitly in the above list).
In fact, it is easy to see that in this case the second rightmost
(resp.\ second highest) path must start with a vertical
(resp.\ horizontal) step, the third rightmost (resp.\ third highest)
path with two vertical (resp.\ horizontal) steps, and the leftmost
(resp.\ lowest) path with $r-1$ vertical (resp.\ horizontal) steps.
Explicitly, we have
\begin{equation}
\det_{1\le i,j,\le r}\big(w(\mathcal P((l+j,k)\to(n_i,m_i)))\big)=
\det_{1\le i,j,\le r}\big(w(\mathcal P((l+j,k+r-j)\to(n_i,m_i)))\big),
\end{equation}
and
\begin{multline}
\det_{1\le i,j,\le r}\big(w(\mathcal P((l,k-j)\to(n_i,m_i)))\big)\\=
\prod_{1\le i<j\le r}w(l+i,k-j)\,
\det_{1\le i,j,\le r}\big(w(\mathcal P((l+j-1,k-j)\to(n_i,m_i)))\big).
\end{multline}
An analogous fact holds if one considers the {\em end} points
instead of the starting points to be consecutive on a horizontal
(resp.\ vertical) line. 
\end{remark}

\section{Identities for multiple elliptic hypergeometric series}\label{secms}

It is straightforward to extend the convolution formulae in
\eqref{ft1}, \eqref{ft2}, and \eqref{ft3}, to the multivariate
setting using the interpretation of nonintersecting lattice paths.
We have the following identities:
\begin{proposition}\label{prop2}
Let $l,k,n,m$ be integers such that $n-l+m-k\ge 0$.

(a) Fix an integer $\nu$ such that $l+r+1\le\nu\le n+1$. Then we have
\begin{multline}\label{ft1m}
\det_{1\le i,j,\le r}\big(w(\mathcal P((l+j,k-j)\to(n+i,m-i)))\big)\\
=\sum_{\underset{t_1\le m-1,t_r\ge k-r}{t_1>t_2>\dots>t_r}}
\det_{1\le i,j,\le r}\big(w(\mathcal P((l+j,k-j)\to(\nu-1,t_i)))\big)
\prod_{s=1}^rw(\nu,t_s)\\\times
\det_{1\le i,j,\le r}\big(w(\mathcal P((\nu,t_j)\to(n+i,m-i)))\big).
\end{multline}

(b) Fix an integer $\nu$ such that $k\le\nu\le m-r$. Then we have
\begin{multline}\label{ft2m}
\det_{1\le i,j,\le r}\big(w(\mathcal P((l+j,k-j)\to(n+i,m-i)))\big)\\
=\sum_{\underset{t_1\ge l+1,t_r\le n+r}{t_1<t_2<\dots<t_r}}
\det_{1\le i,j,\le r}\big(w(\mathcal P((l+j,k-j)\to(t_i,\nu-1)))\big)\\\times
\det_{1\le i,j,\le r}\big(w(\mathcal P((t_j,\nu)\to(n+i,m-i)))\big).
\end{multline}

(c) Fix an integer $\nu$ such that $l+k\le\nu\le n+m$. Then we have
\begin{multline}\label{ft3m}
\det_{1\le i,j,\le r}\big(w(\mathcal P((l+j,k-j)\to(n+i,m-i)))\big)\\
=\sum_{\underset{t_1\ge l+1,t_r\le n+r}{t_1<t_2<\dots<t_r}}
\det_{1\le i,j,\le r}\big(w(\mathcal P((l+j,k-j)\to(t_i,\nu-t_i)))\big)\\\times
\det_{1\le i,j,\le r}\big(w(\mathcal P((t_j,\nu-t_j)\to(n+i,m-i)))\big).
\end{multline}
\end{proposition}

We could also have formulated more general versions
of convolutions where the respective starting and/or end points
of the total paths are not consecutive on antidiagonals
(in the above cases these points are $(l+i,k-i)$ and $(n+i,m-i)$,
$i=1,\dots,r$). However, the advantage of our specific choice is that
all the determinants involved in Proposition~\ref{prop2}
factor into closed form,
by virtue of the determinant evaluations in Proposition~\ref{prop}.
We thus obtain, writing out the identities \eqref{ft1m}, \eqref{ft2m},
and \eqref{ft3m} explicitly, summations which are particularly
attractive since both the summands and the product sides
are completely factored. Each of the above three cases
leads, after suitable substitution of variables, simplification,
and analytic continuation, to the same result. It is a special
case of a multivariate ${}_{10}V_9$ summation formula conjectured
by Warnaar (let $x=q$ in \cite[Cor.~6.2]{Wa}) which has subsequently
been proved by Rosengren~\cite{Ro}.

\begin{theorem}[A multivariate extension of Frenkel and Turaev's
${}_{10}V_9$ summation formula]\label{thmftsum}
Let $a,b,c,d$ be indeterminates, let $m$ be a nonnegative integer,
and $r\ge 1$. Then we have
\begin{multline}\label{nmftsum}
\sum_{0\le k_1<k_2<\dots<k_r\le m}
q^{\sum_{i=1}^r(2i-1)\la_i}\prod_{1\le i<j\le r}
\ta(q^{k_i-k_j},aq^{k_i+k_j})^2\\\times
\prod_{i=1}^r\frac{\ta(aq^{2k_i};p)(a,b,c,d,
a^2q^{3-2r+m}/bcd,q^{-m};q,p)_{k_i}}
{\ta(a;p)(q,aq/b,aq/c,aq/d,bcdq^{2r-2-m}/a,aq^{1+m};q,p)_{k_i}}\\
=q^{-4\binom r3}\left(\frac a{bcdq}\right)^{\binom r2}
\prod_{i=1}^r(q,b,c,d,a^2q^{3-2r+m}/bcd;q,p)_{i-1}\\\times
\prod_{i=1}^r\frac{(q,aq;q,p)_m
(aq^{2-i}/bc,aq^{2-i}/bd,aq^{2-i}/cd;q,p)_{m+1-r}}
{(q,aq/b,aq/c,aq/d,aq^{2-2r+i}/bcd;q,p)_{m+1-i}}.
\end{multline}
\end{theorem}
Note that the Vandermonde determinant-like factor appearing
in the summand of \eqref{nmftsum} is squared.
This distinctive feature is reminiscent of certain Schur function
and multiple $q$-series identities with similar property
(which can also be proved by
the machinery of nonintersecting lattice paths), see e.g.\
\cite[Thms.~5 and 6]{Kr} and \cite[Thms.~27--29]{CEKZ}.

The following result is the natural generalization of 
Theorem~\ref{thmftsum} to the higher level of transformations.
It is a special case of a multivariate ${}_{12}V_{11}$
transformation formula conjectured by Warnaar (let $x=q$ in
\cite[Conj.~6.1]{Wa}) which has subsequently been
proved (in more generality) by Rains~\cite{Ra} and, independently,
by Coskun and Gustafson~\cite{CG}.

\begin{theorem}[A multivariate extension of Frenkel and Turaev's
${}_{12}V_{11}$ transformation formula]\label{thmfttf}
Let $a,b,c,d,e,f$ be indeterminates, let $m$ be a nonnegative integer,
and $r\ge 1$. Then we have
\begin{multline}\label{nmfttf}
\sum_{0\le k_1<k_2<\dots<k_r\le m}
q^{\sum_{i=1}^r(2i-1)k_i}\prod_{1\le i<j\le r}
\ta(q^{k_i-k_j},aq^{k_i+k_j})^2\\\times
\prod_{i=1}^r\frac{\ta(aq^{2k_i};p)(a,b,c,d,e,f,
\la aq^{2-r+m}/ef,q^{-m};q,p)_{k_i}}
{\ta(a;p)(q,aq/b,aq/c,aq/d,aq/e,aq/f,efq^{r-1-m}/\la,aq^{1+m};q,p)_{k_i}}\\
=\prod_{i=1}^r\frac{(b,c,d,ef/a;q,p)_{i-1}}
{(\la b/a,\la c/a,\la d/a,ef/\la;q,p)_{i-1}}\\\times
\prod_{i=1}^r\frac{(aq;q,p)_m\,(aq/ef;q,p)_{m+1-r}\,
(\la q/e,\la q/f;q,p)_{m+1-i}}
{(\la q;q,p)_m\,(\la q/ef;q,p)_{m+1-r}\,
(aq/e,aq/f;q,p)_{m+1-i}}\\\times
\sum_{0\le k_1<k_2<\dots<k_r\le m}
q^{\sum_{i=1}^r(2i-1)k_i}\prod_{1\le i<j\le r}
\ta(q^{k_i-k_j},\la q^{k_i+k_j})^2\\\times
\prod_{i=1}^r\frac{\ta(\la q^{2k_i};p)(\la,\la b/a,\la c/a,\la d/a,e,f,
\la aq^{2-r+m}/ef,q^{-m};q,p)_{k_i}}
{\ta(\la;p)(q,aq/b,aq/c,aq/d,\la q/e,\la q/f,efq^{r-1-m}/\la,
\la q^{1+m};q,p)_{k_i}},
\end{multline}
where $\la=a^2q^{2-r}/bcd$.
\end{theorem}
The $r=1$ case of Theorem~\ref{thmfttf} is
Frenkel and Turaev's ${}_{12}V_{11}$ transformation theorem~\cite{FT},
an elliptic extension of Bailey's ${}_{10}\phi_9$
transformation~\cite[Eq.~(2.9.1)]{GR}.
Again, the Vandermonde determinant-like factor appearing
in the summand of \eqref{nmfttf} is squared. (Similar identities but
with a simple Vandermonde determinant-like factor appearing
in the summand have been derived in \cite{RS1}.)
Due to symmetry the range of summations on both sides of
\eqref{nmfttf} can also be taken over all integers $0\le k_1,\dots,k_r\le m$.
If we let $c=aq/b$ in \eqref{nmfttf}, the left-hand side
reduces to a multivariate $_{10}V_9$ series. On the right-hand side,
since $\la d/a=q^{1-r}$, the sum boils down to just a single term, with
the indices $k_i=i-1$, $1\le i\le r$. The result, after simplifications,
is of course Theorem~\ref{thmftsum}.

It would be particularly interesting to find a combinatorial proof
of \eqref{nmfttf} involving nonintersecting lattice paths.
Even for $r=1$ we so far failed to find a lattice path proof.
We leave this as an open problem.

\section{Future perspectives}\label{secfp}

\subsection{Tableaux and plane partitions}
It is quite clear how one can enumerate objects such as
tableaux or (various classes of) plane partitions with respect to
elliptic weights. First, one has to translate the respective
combinatorial objects via a standard bijection into a set of
nonintersecting lattice paths (see \cite{GV} or \cite{St}).
The translation back, in order to obtain an explicit definition for the
weight of the corresponding combinatorial object, is not difficult. 
In the simplest cases the elliptic generating function is
then expressed, by Theorem~\ref{thnlp}, as a determinant
which may be computed by Proposition~\ref{prop}.
If the starting and/or end points of the lattice paths are
not fixed, one applies instead of Theorem~\ref{thnlp} 
a result by Okada~\cite{O} (see also Stembridge~\cite{St}),
which expresses the generating function as a Pfaffian.
Since the square of a Pfaffian is a determinant of a
skew symmetric matrix, this again involves the computation of
a determinant. It needs to be explored which of the classical
results can be extended to the elliptic setting. Some elliptic
determinant evaluations, other than Warnaar's in Lemma~\ref{lemdet1},
which might be useful in this context have been provided
by Rosengren and present author~\cite{RS2}.

\subsection{Elliptic Schur functions}
One can replace \eqref{ellw} by the more general weight
\begin{equation}\label{ellgw}
w(x;n,m):=
\frac{\ta(ax_m^2q^n,bq^{2n},bq^{2n-1},aq^{1-n}/b,aq^{-n}/b)}
{\ta(aq^n,bx_mq^{2n},bx_mq^{2n-1},ax_mq^{1-n}/b,ax_mq^{-n}/b)}x_m
\end{equation}
(defined on horizontal steps $(n-1,m)\to(n,m)$ of $\mathbb Z^2$),
and enumerate nonintersecting lattice paths, corresponding to tableaux,
with respect to \eqref{ellgw}. The result is an elliptic extension of
Schur functions (which perhaps are no longer orthogonal with respect
to any elliptic scalar product) which, when ``principally specialized''
($x_i\mapsto q^i$, $i\ge 0$) by construction
factors into closed form in view of Proposition~\ref{prop}.
It should be worth investigating whether these elliptic Schur functions
have other nice properties (as they do have in the classical
case, see \cite{M}). As a matter of fact, they do not seem to be
related to (the $t=q$ cases of) any of the $BC$-symmetric functions
considered in \cite{CG} or \cite{Ra}. On the other hand, it would be
already interesting to study limiting cases of the $p=0$ case of
these elliptic Schur functions. One would hope that the
Hall--Littlewood functions (which are an important one parameter
extension of the Schur functions, cf.~\cite{M})
would then appear as a special case, which would then admit
a surprising combinatorial interpretation in terms of
lattice paths. Unfortunately, as a matter of fact,
the Hall--Littlewood functions do not seem to be contained
in the above considered family of elliptic Schur functions.

\subsection{Other weight functions}
We were able to disguise Frenkel and Turaev's ${}_{10}V_9$ summation
formula as a convolution identity of elliptic binomial coefficients
(see also Rains~\cite[Sec.~4]{Ra} and Coskun and Gustafson~\cite{CG}).
In our case this involved lattice paths with respect to elliptic weights.
Similarly, it should also be feasible to reproduce other
known convolution formulae (such as Abel's generalization of the
binomial theorem or the Hagen--Rothe summation, cf.\ \cite{S},
or others) using lattice paths with appropriately chosen weights.
The three types of convolutions, displayed in \eqref{ft1},
\eqref{ft2}, and \eqref{ft3}, still hold, but may then lead to
mutually different identities. One can also try to work with
{\em bibasic} weights (either elliptic or non-elliptic),
in order to recover some of the identities in \cite[Secs.~3.6 and 3.8]{GR}
and in \cite{Wa}. It seems likely that in the non-elliptic case
(here we mean that there is no nome $p$, or $p=0$)
Bill Gosper used exactly this method to first derive his
``strange evaluations'' (which were later subsumed/generalized in
\cite[Secs.~3.6 and 3.8]{GR}).
Of course, whatever identities or other results one obtains
by lattice path interpretation, one can check for possible related
determinant evaluations. Also the other direction should be
investigated, e.g.\ does Warnaar's quadratic elliptic determinant
in \cite[Thm.~4.17]{Wa} correspond to a specific set of nonintersecting
lattice paths with quadratic elliptic weight function?

\subsection{``Elliptic'' combinatorics}
I strongly believe that the results presented in this paper do not stand
alone, i.e., that elliptic enumeration is not necessarily restricted to
lattice paths.
In the same way as the generating functions for various classes of
combinatorial objects, most notably, of partitions, which correspond
to paths, can be expressed in terms of $q$-series,
closed form elliptic generating functions for several of these classes
should exist as well.
The main idea would be to replace $q$-weights by suitable
elliptic weights, and then try to make the further analysis work out.
There are certainly restrictions to the elliptic approach
(besides that the objects counted should be finite).
For instance, still considering paths in $\mathbb Z^2$,
Andr\'e's reflection principle (cf.\ \cite[p.~22]{C})
is not applicable as it is not anymore weight invariant.
Techniques involving shifting paths (as in \cite[Prop.~1]{GK}),
however, may still work with delicate handling (see Remark~\ref{remsh}).
Besides lattice path enumeration, a good area where to look for
elliptic extensions would presumably be a general combinatorial theory
such as Viennot's theory of heaps~\cite{V}.

\bibliographystyle{amsalpha}

\end{document}